\documentclass[a4paper,11pt,english]{smfart}

\usepackage[utf8]{inputenc}

\usepackage{amsthm}
\usepackage{amsmath,amssymb,amsfonts}
\usepackage{epsfig,color}
\usepackage{mathrsfs}
\usepackage{stmaryrd}
\usepackage{hyperref}

\newcommand{\LL}{\textbf{L}}
\newcommand{\I}{\textbf{Id}}
\newcommand{\E}{\mathbb{E}}
\newcommand{\R}{\mathbb{R}}

\newcommand{\s}{\scalebox{1.3}{${}^\sharp$ \hspace{-0.15 cm}}}
\newtheorem{thm}{Theorem}

\title{Regularization along central convergence on second and third Wiener chaoses}

\author{Guillaume Poly} 
\address{IRMAR, Universit\'e de Rennes 1} 
\email{guillaume.poly@univ-rennes1.fr} 
\urladdr{https://sites.google.com/site/guillaumejpoly/home}

\usepackage{fancyhdr}
\renewcommand{\headheight}{13pt}

\begin{document}

\begin{abstract}
Consider $F$ an element of the second Wiener chaos with variance one. In full generality, we show that, for every integer $p\ge 1$, there exists $\eta_p>0$ such that if $\kappa_4(F)<\eta_p$ then the Malliavin derivative of $F$ admits a negative moment of order $p$. This entails that any sequence of random variables in the second Wiener chaos converging in distribution to a non--degenerated Gaussian is getting more regular as its distribution is getting close to the normal law. This substantially generalizes some recent findings contained in \cite{hu2014convergence,hu2015density,nourdin2016fisher} where analogous statements were given with additional assumptions which we are able to remove here. Moreover, we provide a multivariate version of this Theorem.

\smallskip

Our main contribution concerns the case of the third Wiener chaos which is notoriously more delicate as one cannot anymore decompose the random variables into a linear combination of i.i.d. random variables. We still prove that the same phenomenon of regularization along central convergence occurs. Unfortunately, we are not able to provide a statement as strong as the previous one, but we can show that the usual non--degeneracy estimates of the Malliavin derivative given by the Carbery-Wright inequality can be improved by a factor three. Our proof introduces new techniques such that a specific Malliavin gradient enabling us to encode the distribution of the Malliavin derivative by the spectrum of some Gaussian matrix. This allows us to revisit the fourth moment phenomenon in terms of the behavior of its spectral radius.
\end{abstract}

\maketitle

\section{The case of the second Wiener chaos}

Let us consider $F$ an element of the second Wiener chaos $\text{Ker}(\LL+2\I)$ where $\LL$ stands for the Ornstein-Uhlenbeck operator. We refer the reader to the book \cite{bakry2013analysis} for an introduction to the formalism of Wiener chaoses from the point of view of Markov generators. It is a standard fact that one may find $\{\alpha_{k}\}_{k\ge 1}$ a sequence of real numbers in $\textit{l}^2(\mathbb{N}^\star)$ such that

\begin{equation}\label{sec.chaos.decomp}
F=\sum_{k=1}^\infty \alpha_{k} \left(G_k^2-1\right).
\end{equation}

Above, $\{G_k\}_{k\ge 1}$ stands for an \textbf{i.i.d.} sequence of standard Gaussian random variables which belongs to the first Wiener chaos. The reader may for instance consult \cite[2.7.13]{nourdin2012normal} for a corresponding proof. We shall assume that $\E\left(F^2\right)=1$, it is then a well known fact that $\kappa_4(F)=\E(F^4)-3$ controls the proximity in distribution between $F$ and the standard Gaussian law. Our main goal is to show that, provided that $\kappa_4(F)$ is small enough, the regularity of the distribution $f$ of $F$ increases. The key ingredient is to establish negative moments for the Malliavin derivative of $F$, namely

\begin{equation}\label{negativemoments}
\forall p\ge 1, \exists \delta_p>0\\,\,\text{s.t.}\,\,\, \kappa_4(F)<\delta_p\,\,\Rightarrow \frac{1}{\Gamma[F,F]}\in L^p(\Omega,\mathcal{F},\mathbb{P}).
\end{equation}

Actually we will prove a more quantitative version of the previous statement which is given in the Theorem below.

\begin{thm}\label{Mainsecond}
Let $F\in\text{Ker}(\LL+2\I)$ which satisfies $\E(F^2)=1$. One has, for every $p\ge 1$,
\begin{equation}\label{negativemomentsbis}
\kappa_4(F)<\frac{24}{2^p (p+1)!}\,\Rightarrow\,\frac{1}{\Gamma[F,F]}\in \bigcap_{q<\frac{p}{2}} L^q(\Omega,\mathcal{F},\mathbb{P}).
\end{equation}
\end{thm}

The latter considerably improves recent findings on this topic for the case of the second Wiener chaos since it completely removes any kind of additional assumptions and only requires the central convergence since the integrability of the inverse of the Malliavin derivative is directly related to the smallness of the fourth cumulant. One can for instance consult \cite[Thm 7.2 \text{and} Thm 7.3]{hu2014convergence}, \cite[Thm 1.5 : case of the second Wiener chaos]{hu2015density}, \cite{nourdin2016fisher} or else the nice survey \cite[page 377]{crisan2014stochastic} for related statements, all of them requiring additional assumptions ensuring negative moments for the Malliavin derivative.

\medskip

\begin{proof} First of all, using the representation (\ref{sec.chaos.decomp}), we immediately deduce that

$$\Gamma[F,F]=4\sum_{k=1}^\infty \alpha_k^2 G_k^2;$$

Computing the Laplace transform of the previous expression gives

\begin{equation}\label{Lapl.}
\forall \lambda>0,\,\E\left(\exp\left(-\lambda \Gamma[F,F]\right)\right)=\prod_{k=1}^\infty \frac{1}{\sqrt{1+8 \lambda \alpha_k^2}}.
\end{equation}

Now, let us introduce the symmetric elementary functions
\begin{equation}\label{symfunc}
S_{p}=\sum_{i_1<i_2<\cdots<i_p} \alpha_{i_1}^2\alpha_{i_2}^2\cdots \alpha_{i_p}^2,
\end{equation}
as well as the so-called \textit{Newton sums} which are in turn expressed in terms of cumulants:

\begin{equation}\label{Newt-Cum}
\mathcal{N}_p=\sum_{k=1}^\infty \alpha_k^{2p}=\frac{1}{2^{2p-1}(2p-1)!} \kappa_{2p}(F).
\end{equation}

It is a standard fact that symmetric elementary functions can be expressed in terms of the Newton sums. In particular, one has the following identities, which are due to Newton and Girard and whose proof may be found for instance in \cite{richter1949einfacher}.

\begin{equation}\label{Newt-Sym-rel}
S_p=(-1)^p \sum_{m_1+2m_2+\cdots+pm_p=p} \prod_{i=1}^p \frac{(-\mathcal{N}_i)^{m_i}}{m_i! i^{m_i}}
\end{equation}

One needs to estimate the denominator of the right hand side of equation (\ref{Lapl.}). A simple expansion gives
\begin{eqnarray*}
\prod_{k=1}^\infty (1+8\lambda \alpha_k^2)&=&1+\sum_{p=1}^\infty 8^p \lambda^p S_p.\\
\end{eqnarray*}
Now we will isolate in the formula (\ref{Newt-Sym-rel}) the case $m_1=p$, which forces $m_2=m_3=\cdots=m_p=0$.
\begin{eqnarray*}
S_p&=&\frac{\mathcal{N}_1^p}{p!}+(-1)^p \sum_{\begin{subarray}{c}
m_1<p\\m_1+2m_2+\cdots+pm_p=p
\end{subarray}} \prod_{i=1}^p \frac{(-\mathcal{N}_i)^{m_i}}{m_i! i^{m_i}}\\
&=&\frac{1}{2^p p!}+(-1)^p \sum_{\begin{subarray}{c}
m_1<p\\m_1+2m_2+\cdots+pm_p=p
\end{subarray}} \prod_{i=1}^p \frac{(-\mathcal{N}_i)^{m_i}}{m_i! i^{m_i}}.\\
\end{eqnarray*}

We notice that the condition $m_1<p$ implies that there exists $i\in\{2,\cdots,p\}$ such that $m_i>0$. Besides, since for any $i\ge 1$ we have $\alpha_i^2\le 1$ we get that for all $i\ge 2$ we have $\mathcal{N}_i\le \mathcal{N}_2=\frac{1}{48}\kappa_4(F)$ as well as $\mathcal{N}_i\le 1$ for every $i\ge 2$. We have then the rough inequalities
\begin{eqnarray*}
\left|\sum_{\begin{subarray}{c}
m_1<p\\m_1+2m_2+\cdots+pm_p=p
\end{subarray}} \prod_{i=1}^p \frac{(-\mathcal{N}_i)^{m_i}}{m_i! i^{m_i}}\right|&\le& \frac{\kappa_4(F)}{48}\sum_{m_1+2m_2+\cdots+pm_p=p} \prod_{i=1}^p \frac{1}{m_i! i^{m_i}}\\
&\le& p \frac{\kappa_4(F)}{48}
\end{eqnarray*}

The last inequality above comes from the fact that
\begin{eqnarray*}
\sum_{m_1+2m_2+\cdots+pm_p=p} \prod_{i=1}^p \frac{1}{m_i! i^{m_i}}&\le& \prod_{i=1}^p \sum_{m_i=0}^p \frac{1}{m_i! i^{m_i}}\\
&\le& \prod_{i=1}^p e^{\frac{1}{i}}\\
&\le& p
\end{eqnarray*}

There is certainly room for improvement in the aforementioned inequalities but we are not focused here on establishing sharp bounds. As a result, we have established for all $p\ge 1$ the following inequality
\begin{equation}\label{ineq-newt}
\left|S_p-\frac{1}{2^p p!}\right|\le \frac{p}{48} \kappa_4(F).
\end{equation}
Hence, as soon as $\frac{p}{48}\kappa_4(F)\le \frac{1}{p!}\left(\frac{1}{2^p}-\frac{1}{2^{p+1}}\right)$ one has $S_p\ge \frac{1}{2^{p+1}p!}$. This holds when $\kappa_4(F)<\frac{24}{2^p (p+1)!}$ and it implies
\begin{equation*}
\forall \lambda>0,\,\,\prod_{k=1}^\infty (1+8\lambda \alpha_k^2)\ge 1+\frac{1}{2^{p+1}p!} 8^p\lambda^p=1+\frac{4^p}{2~p!} \lambda^p.
\end{equation*}
and thus
\begin{equation}
\forall \lambda>0, \E\left(\exp\left(-\lambda \Gamma[F,F]\right)\right)=\frac{1}{\sqrt{1+\frac{1}{2~p!}4^p \lambda^p}}\le \frac{\sqrt{2~p!}}{2^p \lambda^{\frac{p}{2}}}.
\end{equation}
Now, take $\epsilon>0$, one has
\begin{eqnarray*}
\mathbb{P}\left(\Gamma[F,F]<\epsilon\right)&=&\mathbb{P}\left(\exp\left(-\lambda \Gamma[F,F]\right)\ge \exp(-\lambda \epsilon)\right)\\
&\le& e^{\lambda \epsilon} \E\left(e^{-\lambda \Gamma[F,F]}\right)\\
&\le& \frac{\sqrt{2~p!}}{2^p} \frac{e^{\lambda \epsilon}}{\lambda^{\frac{p}{2}}}.
\end{eqnarray*}
Setting $\lambda=\frac{1}{\epsilon}$ gives that 
$$\mathbb{P}\left(\Gamma[F,F]<\epsilon\right)\le \frac{\sqrt{2~p!}}{2^p} \epsilon^{\frac{p}{2}},$$
and ensures that $1/\Gamma[F,F]\in L^q$ for every $q<\frac{p}{2}$ which achieves the proof.
\end{proof}

\section{Multivariate case in the second Wiener chaos}

In this section, we extend to the multivariate case the content of Theorem \ref{Mainsecond}. We will write it in a sequential way since it is more convenient for us, though it dos not provide a quantitative statement. Consider $F_n=(F_{n,1},\cdots,F_{n,d})$ such that for every $i\in\{1,\cdots,d\}$ we have $F_{n,i} \in \text{Ker}(\LL+2 \I)$. Moreover we will also assume that $\text{Cov}(F_n)=\I_d$. Then we shall prove the following Theorem:

\begin{thm}\label{multi-second}
Assume that for every $i\in\{1,\cdots,d\}$ we have $\kappa_4(F_{n,i})\to 0$. Then, for every $q\ge 1$ there exists $C_q>0$ and $n_q\ge 1$ such that

\begin{equation}
\E\left( e^{i \sum_{k=1}^d \xi_k F_{k,n} } \right) \le \frac{C_q}{\left(\xi_1^2+\cdots+\xi_d^2\right)^{\frac q 2}}.
\end{equation}
\end{thm}

Using classical Fourier inversion methods, this shows that the sequence of the joint densities of $(F_{n,1},\cdots,F_{n,d})$ converges towards the Gaussian density in $\mathcal{C}^\infty$ topology.

\begin{proof}
Usually, to handle the non--degeneracy of a random vector by Malliavin calculus, one needs to show negative moments for the determinant of its Malliavin matrix. Here, we were not able to prove that $\det\left(\left(\Gamma[F_{n,i},F_{n,j}\right)_{1\le i,j\le d}\right)$ admits negative moments for $n$ large enough. Instead of that we will only work somehow with univariate variables.

\medskip

Using for instance Peccati-Tudor Theorem \cite{peccati2005gaussian} , whose an alternate proof may be found in \cite{campese2016multivariate}, we know that for every $i\ne j$, $\Gamma[F_{n,i},F_{n,j}]\to 0$ in $L^2(\mathbb{P})$. As a result, for every $(t_1,\cdots,t_d)\in \mathcal{S}^{d-1}$ the unit sphere of $\R^d$ we have
\begin{eqnarray*}
&&\text{Var}\left(\Gamma\left[\sum_{i=1}^d t_i F_{n,i},\sum_{i=1}^d t_i F_{n,i}\right]\right)= \text{Var}\left(\sum_{i,j=1}^d t_i t_j \Gamma\left[F_{n,i},F_{n,j}\right]\right)\\
&=& \sum_{i=1}^d t_i^4 \text{Var}\left(\Gamma[F_{n,i},F_{n,i}]\right)+\sum_{i\neq j \,\text{or}\, k\neq l} t_i t_j t_k t_l \text{Cov}\left(\Gamma[F_{n,i},F_{n,j}],\Gamma[F_{n,k},F_{n,l}]\right)\\
\end{eqnarray*}
As a result we get
\begin{equation}\label{control1multi}
\text{Var}\left(\Gamma\left[\sum_{i=1}^d t_i F_{n,i},\sum_{i=1}^d t_i F_{n,i}\right]\right)\le \max_{1\le i\le d} \text{Var}\left(\Gamma[F_{n,i},F_{n,i}]\right)+d^2 \max_{i\neq j}\|\Gamma[F_{n,i},F_{n,j}]\|_2.
\end{equation}
As the right hand side is independent of $(t_1,\cdots,t_d)$ and since on any Wiener chaos $\text{Var}(\Gamma[F_n,F_n])\to 0$ and $\kappa_4(F_n)\to 0$ are equivalent statements, one gets that

\begin{equation}\label{control2multi}
\max_{(t_1,\cdots,t_d)\in\mathcal{S}^{d-1}}\kappa_4\left(\sum_{i=1}^d t_i F_{n,i}\right)\to 0.
\end{equation}

Thus, in virtue of Theorem \ref{Mainsecond}, for any $q\ge 1$ there exists $n_q\ge 1$ such that 
$$\forall n\ge n_q,\,\forall (t_1,\cdots,t_d) \in \mathcal{S}^{d-1},\,\frac{1}{\Gamma\left[\sum_{i=1}^d t_i F_{n,i},\sum_{i=1}^d t_i F_{n,i}\right]} \in L^p\left(\Omega,\mathcal{F},\mathbb{P}\right).$$

Now, let us fix $\textbf{t}=(t_1,\cdots,t_d)\in\mathcal{S}^{d-1}$ and let us introduce for convenience the notation $F_n(\textbf{t})=\sum_{i=1}^d t_i F_{n,i}.$ Then, using integrations by parts techniques (see for instance \cite[Prop 2.1.4]{nualart2006malliavin}) and combining the uniformity in the bound (\ref{control2multi}) with Theorem \ref{Mainsecond}, one can show that there is a constant $C_p>0$ such that for every $\phi\in\mathcal{C}^{p-1}_b(\R)$: 
$$\forall \textbf{t}\in\mathcal{S}^{d-1},\,\forall n\ge n_p,\,\left|\E\left(\phi^{p-1}(F_n(\textbf{t})\right)\right|\le C_p \|\phi\|_\infty.$$

Applying this to $\textbf{t}_\xi=\frac{1}{\sqrt{\xi_1^2+\cdots+\xi_d^2}}(\xi_1,\cdots,\xi_d)$ and $\phi(x)=e^{i x \sqrt{\xi_1^2+\cdots+\xi_d^2}}$ gives that

\begin{eqnarray*}
\left|\E\left( e^{i \sum_{k=1}^d \xi_k F_{k,n}}\right)\right|&=&\left|\E\left( \phi( F_n(\textbf{t}_\xi) \sqrt{\xi_1^2+\cdots+\xi_d^2}\right)\right|\\
\le \frac{C_p}{\left(\xi_1^2+\cdots+\xi_d^2\right)^{\frac{p-1}{2}}}.
\end{eqnarray*}
\end{proof}

\section{The case of the third Wiener chaos}

Let us consider here $F$ an element of $\text{Ker}(\LL+3\I)$ which satisfies $\E(F^2)=1$. Contrarily to the case of the second Wiener chaos, it is not anymore possible to decompose $F$ as a linear combination of independent random variables. This issue remains for higher order chaoses and mainly explains why it is in general very difficult to extend to these settings some results holding true for the second Wiener chaos. For instance, one still ignores what are the possible limits in distribution of chaotic random variables or order strictly greater than $2$. In contrast, we know that every limit in distribution of random variables in the second Wiener chaos may be written as an element of the second Wiener chaos plus some independent Gaussian random variable, see for instance \cite[Thm 3.1]{nourdin2012convergence} or else \cite[Thm 1.2]{bogachev2015two} for the multivariate case. Another important instance where the second Wiener chaos provides stronger results is given by \cite[Thm 3.4]{nourdin2012convergence}, where for a wide class of targets, the convergence in distribution is shown to be equivalent to the convergence of a finite number of moments/cumulants. For sequences lying in Wiener chaoses of order larger than $3$, when the target is not Gaussian or else Gamma, we still don 't know whether the convergence in distribution is ensured by the convergence of a finite number of moments/cumulants. Closely related to these topics, one might read \cite{bourguin2018four} where quantitative fourth moment theorems are given for targets lying in the Pearson family or else \cite{dobler2018fourth} for analogous questions in the Poisson space.

\smallskip

In this section, we investigate whether the phenomenon of regularization along central convergence,  which is established for the second Wiener chaos in the last section, still holds for the third chaos. Unfortunately, we are not able to provide a statement as strong as Theorem \ref{Mainsecond}. Nevertheless, we shall prove that provided that $\kappa_4(F)$ is small enough, one can improve by a factor tending to three the usual estimates of the non degeneracy of $\Gamma[F,F]$  given by the Carbery-Wright inequality. Concretely, we will establish that

\begin{thm}
For every $\theta < \frac{3}{4}$, there exits $\epsilon>0$ such that
$$\kappa_4(F)<\epsilon \Rightarrow \frac{1}{\Gamma[F,F]^\theta} \in L^1(\Omega,\mathcal{F},\mathbb{P}).$$
\end{thm}

This should be compared with the threshold $\frac 1 4$ which is given by a standard application of Carbery-Wright inequality. It should also be emphasized that we will not provide a quantitative proof of the previous result, since it seemed to us too technical and beyond the scope of this article. Abstract constants proceed from a reasoning by the absurd at the end of the third step of the proof below. Finally let us mention that Carbery-Wright inequality is a key ingredient in a serie of recent papers as it enables one to provide quantitative statements in terms of total variation, see for instance the following highly non exhaustive list \cite{bally2017total,nourdin2013convergence,bogachev2018fractional}. As such, any improvement of the Carbery-Wright result in many applications regarding the current state of this art.

\medskip

\begin{proof}

\medskip

\underline{Step 0: preliminary material and notations}

\medskip

\smallskip

Let us consider $(X_i)_{i\ge 1}$ an i.i.d. sequence of standard Gaussian random variables and the underlying probability space $(\Omega,\mathcal{F},\mathbb{P})$. Without loss of generality, we shall assume that $\mathcal{F}=\sigma(X_i;i\ge 1)$. Let us briefly recall that the Malliavin operators on the Wiener space $(\Omega,\mathcal{F},\mathbb{P})$ are defined in the following way:

\begin{eqnarray*}
&&\forall f \in \mathcal{C}^1_b(\R^d,\R),\\
&&\Gamma[f(X_1,\cdots,X_d),f(X_1,\cdots,X_d)]=\sum_{1\le i,j \le d} \partial_i f(X_1,\cdots,X_d)^2\\
&&\forall f \in \mathcal{C}^2_b(\R^d,\R),\\
&& L\left[f(X_1,\cdots,X_d)\right]=\Delta f (X_1,\cdots,X_d) -\sum_{i=1}^d X_i \partial_i f (X_1,\cdots,X_d).
\end{eqnarray*}

Consider $\{a(i,j,k)\}_{1\le i,j,k\le N}$ a sequence of real numbers indexed by $\llbracket 1, N \rrbracket^3$ which is assumed to be symmetric: $\forall \sigma \in \mathcal{S}_3,$ and every $(i_1,i_2,i_3) \in \llbracket 1, N \rrbracket^3$ one has $a(i_{\sigma(1)},i_{\sigma(2)},i_{\sigma(3)})=a(i_1,i_2,i_3)$. We will also assume that whenever $i_1=i_2$ or $i_2=i_3$ or else $i_1=i_3$ we have $a(i_1,i_2,i_3)=0$. Then, let us introduce
\begin{eqnarray*}
F&:=&F(X_1,\cdots,X_N)\\
&=&\sum_{1\le i_1,i_2,i_3 \le N} a(i_1,i_2,i_3) X_{i_1} X_{i_2} X_{i_3} \in \textbf{Ker}(\LL+3 \I).
\end{eqnarray*}

We will further assume that

$$\E(F^2)=6 \sum_{1\le i_1 < i_2<i_3\le N} a(i_1,i_2,i_3)^2=1.$$

In the sequel we will work with $F$ which is a homogeneous polynomial of degree $3$ into a finite numbers of Gaussian random variables. All the forthcoming estimates will be independent on the number of entries $N$, and therefore can be extended verbatim to the infinite dimensional setting.

\medskip

Finally the proof requires here a particular \textit{gradient}, customarily called the sharp operator, which has been introduced by N. Bouleau, see \cite[page 135]{bouleau2010dirichlet} or \cite[page 80]{bouleau2003error}. This gradient is particularly convenient since it maintains somehow the \textit{Gaussian structure} with opposition with the standard choice of the literature to introduce an auxiliary Hilbert space of the form $L^2([0,T])$. It plays a crucial role in our approach. To do so, we need a copy $(\hat{\Omega},\hat{\mathcal{F}},\hat{\mathbb{P}})$ of $(\Omega,\mathcal{F},\mathbb{P})$ as well as $(\hat{X}_i)_{i\ge 1}$ a corresponding i.i.d. sequence of standard Gaussian such that $\hat{\mathcal{F}}=\sigma(\hat{X}_i;i\ge 1)$.  For any $m\ge 1$ and any $F\in\mathcal{C}^1_{\text{Pol}}(\R^m,\R)$, the set of functions of $\mathcal{C}^1(\R^m,\R)$ whose gradient has a polynomial growth, one may define the following \textit{sharp operator}:
\begin{equation}\label{sharp operator}
\s F(X_1,\cdots,X_m)
:=\sum_{i=1}^m \partial_i F(X_1,\cdots,X_m) \hat{X_i}.
\end{equation}

\medskip

\underline{Step 1) Encoding the law of $\Gamma[F,F]$ by the spectrum of a Gaussian matrix}

\medskip

For the sake of clarity, $\E$ denotes the expectation with respect to $\mathbb{P}$ and $\hat{\E}$ stands for the expectation with respect to $\hat{\mathbb{P}}$. Relying on the definition (\ref{sharp operator}) and the formula for the square field operator which was recalled previously, one deduces that $\Gamma[F,F]=\hat{\E}((\s F)^2)$. Moreover, conditionally to $(X_1,\cdots,X_N)$, the random variable $\s F$ is a centered Gaussian whose variance is $\Gamma[F,F]$. As a result we get the formula

\begin{equation}\label{Lapl-caract-gradient}
\E\left(e^{-\frac{\Gamma[F,F]}{2} \xi^2}\right)=\E\hat{\E}\left(e^{i \xi \s F}\right)
\end{equation}

This relates the Laplace transform of the square field operator with the characteristic function of the sharp gradient. On the other hand one may write

\begin{eqnarray*}
\s F&=& \sum_{1\le i_1,i_2,i_3\le N} a(i_1,i_2,i_3) \s \left(X_{i_1} X_{i_2} X_{i_3} \right)\\
&=& \sum_{1\le i_1,i_2,i_3\le N} a(i_1,i_2,i_3) \left(\s X_{i_1} X_{i_2} X_{i_3}+X_{i_1} \s X_{i_2} X_{i_3}+X_{i_1} X_{i_2} \s X_{i_3} \right)\\
&=&\sum_{1\le i_1,i_2,i_3\le N} a(i_1,i_2,i_3) \left(\hat{X}_{i_1} X_{i_2} X_{i_3}+X_{i_1} \hat{X}_{i_2} X_{i_3}+X_{i_1} X_{i_2} \hat{X}_{i_3} \right)\\
&=& \sum_{1\le i_1,i_2 \le N} \hat{\alpha}(i_1,i_2) X_{i_1} X_{i_2}
\end{eqnarray*}
where, by using symmetry properties of the sequence $\{a(i,j,k)\}$, we have

$$\hat{\alpha}(i_1,i_2)=3 \sum_{1\le k \le N} a(i_1,i_2,k)\hat{X}_{k}.$$

The matrix $\hat{A}\in\mathcal{M}_{N}(\mathbb{R})$ whose entries are $\left(\hat{\alpha}(i_1,i_2)\right)_{1\le i_1,i_2\le N}$ is then a symmetric matrix with Gaussian entries. We can diagonalize it and find $N$ random variables $(\hat{\lambda}_1,\cdots,\hat{\lambda}_N)$ which correspond to the spectrum of $\hat{A}$ and are ordered in the following way: $|\hat{\lambda}_1|\ge|\hat{\lambda}_2|\ge\cdots\ge |\hat{\lambda}_N|$. There also exists $\hat{P}$ a random variable taking values in $ \mathcal{O}_N(\R)$, the set of orthonormal matrix, such that $\hat{A}=\hat{P}\,\text{Diag}(\hat{\lambda_1},\cdots,\hat{\lambda}_N) {}^t \hat{P}$. As a result, denoting by $(Y_1,\cdots,Y_N)=(X_1,\cdots,X_N) P$, we have
\begin{eqnarray*}
\sum_{1\le i_1,i_2 \le N} \hat{\alpha}(i_1,i_2) X_{i_1} X_{i_2}&=& (X_1,\cdots,X_N) \hat{A} {}^t (X_1,\cdots,X_N)\\
&=& (Y_1,\cdots,Y_N) \text{Diag}(\hat{\lambda_1},\cdots,\hat{\lambda}_N) {}^t (Y_1,\cdots,Y_N)\\
&=&\sum_{k=1}^N \hat{\lambda}_k Y_k^2\\
&=&\sum_{k=1}^N\hat{\lambda}_k (Y_k^2-1),
\end{eqnarray*}
where the last equality comes from the fact that the trace of $\hat{A}$ is zero. Let us note that, since $\hat{P}$ is an orthonormal matrix which is independent of $(X_1,\cdots,X_N)$, we get that $(Y_1,\cdots,Y_N)$ is a vector of i.i.d. standard Gaussian conditionally to $(\hat{X}_1,\cdots,\hat{X}_N)$. As a result, in virtue of the Fubini Theorem and the conditional independence aforementioned, one gets
\begin{eqnarray*}
\E\hat{\E} \left( e^{i \xi \s F}\right)&=&\hat{\E}\E\left( e^{i \xi \s F}\right)\\
&=& \hat{\E} \left( \prod_{k=1}^N \E \left(e^{i \xi \hat{\lambda}_k (Y_k^2-1)} \right) \right)\\
&=&\hat{\E}\left( \prod_{k=1}^N \frac{1}{\sqrt{1-2 i \hat{\lambda}_k \xi}}\right)
\end{eqnarray*}
Note that, $\sqrt{1+i x}$ may be defined unambiguously for every $x\in\R$ by $\sqrt{1+ix}=(1+x^2)^{\frac{1}{4}} e^{i \arctan(x)/2}$. Coming back to the equation (\ref{Lapl-caract-gradient}) we get the following relation between the square field operator and the spectrum of $\hat{A}$:
\begin{equation}\label{Gamma-spec}
\forall \xi \in \mathbb{R},\,\,\E\left(e^{-\frac{\Gamma[F,F]}{2} \xi^2}\right)=\hat{\E}\left(\frac{1}{\sqrt{\det\left(\I-2 i \xi \hat{A}\right)}}\right).
\end{equation}
As a matter of fact, the distribution of $\Gamma[F,F]$ is fully encoded by the distribution of the spectrum of $\hat{A}$.

\medskip

\underline{Step 2) Bounding the spectral radius of $\hat{A}$ by the fourth cumulant of $F$:}

\medskip

In this step we use tail estimates of Gaussian chaoses which have been proved for instance in \cite[Thm1 or equation (2)]{latala2006estimates}. We will need them only for the case of the second Wiener chaos. These estimates assert that there exists an absolute constant $C>1$ such that, for every $N\ge1$, any symmetric real matrix $M\in\mathcal{M}_N$ and any $p\ge 1$:

\begin{equation}\label{Gaussian-Tail}
\frac{1}{C} \left(\sqrt{p}\sqrt{\text{Tr}(M^2)}+p\rho_M\right)\le \left\|\sum_{i,j\le N}M_{i,j} X_i X_j\right\|_p\le C \left(\sqrt{p}\sqrt{\text{Tr}(M^2)}+p\rho_M\right)
\end{equation}
\smallskip
where $\rho_M$ stands for the spectral radius of $M$. Applying this in the case $M=\hat{A}$ and considering the $2p$-norm with respect to $\mathbb{P}$ only gives the following inequality $\hat{\mathbb{P}}$-almost surely.
\begin{eqnarray*}
\E\left( \left|\s F\right|^{2p}\right)\ge \left(\frac{4}{C^2}\right)^p |\hat{\lambda}_1|^{2p} p^{2p}.
\end{eqnarray*}
Taking now the expectation with respect to $\hat{\E}$ gives
\begin{eqnarray*}
\hat{\E}\E\left(\left|\s F\right|^{2p}\right)&=&\E\hat{\E}\left(\left|\s F\right|^{2p}\right)\\
&=&\E\left(\Gamma[F,F]^p\right) \frac{(2p)!}{2^p p!}\\
&\ge& \left(\frac{4}{C^2}\right)^p p^{2p} \hat{\E}\left(|\lambda_1|^{2p}\right).
\end{eqnarray*}
The last equality uses the fact that $\s F$ is Gaussian of variance $\Gamma[F,F]$ conditionally to $\hat{\mathcal{F}}$ and $\frac{(2p)!}{2^p p!}$ is the moment of order $2p$ of a standard Gaussian random variable. Note that by Stirling formula one has
\begin{eqnarray*}
\frac{(2p)!}{2^p p!} \frac{1}{p^{2p}}&\sim& \sqrt{2}\frac{2^p}{e^p} p^p \frac{1}{p^{2p}}\\
&\sim&\sqrt{2}\left(\frac{2}{e}\right)^p \frac{1}{p^p}
\end{eqnarray*}
As a result, for some universal constant $C$ we have
\begin{equation}\label{tail-final1}
\left\|\hat{\lambda}_1\right\|_{2p}\le \frac{C}{\sqrt{p}}\sqrt{\left\|\Gamma[F,F]\right\|_p}.
\end{equation}
On the other hand, by the triangle inequality and hypercontractivity on Wiener chaoses, we get

\begin{eqnarray*}
\left\|\Gamma[F,F]\right\|_p&\le& \left\|\Gamma[F,F]-3\right\|_p+3\\
&\le& (p-1)^2 \left\|\Gamma[F,F]-3\right\|_2+3\\
&\le& p^2 \sqrt{\text{Var}\left(\Gamma[F,F]\right)}+3
\end{eqnarray*}
Besides, a central result of Nourdin-Peccati theory, applied here for the third Wiener chaos case, asserts that 
$$ \sqrt{\text{Var}\left(\Gamma[F,F]\right)}\le 3 \sqrt{\kappa_4(F)}.$$

Many proofs may be found in the literature, one can consult for instance \cite[Thm3.2]{azmoodeh2014fourth} for a proof or else \cite{nourdin2012normal} for an overview of the literature around this type of inequalities.
Gathering all theses facts gives, for some universal constant $C$ which may change from line to line, that
\begin{eqnarray*}
\left\|\hat{\lambda}_1\right\|_{2}&\le&\left\|\hat{\lambda}_1\right\|_{2p}\\
&\le& \frac{C}{\sqrt{p}} \sqrt{p^2 \sqrt{\text{Var}\left(\Gamma[F,F]\right)}}+\frac{\sqrt{3 C}}{\sqrt{p}}\\
&\le& C\left(\sqrt{p} \kappa_4(F)^{\frac 1 4}+\frac{1}{\sqrt{p}}\right)
\end{eqnarray*}
Optimizing in $p$ then gives, for some absolute constant $C>0$ that

\begin{equation}\label{finalboundspec}
\left\|\hat{\lambda}_1\right\|_{2}\le C \kappa_4(F)^{\frac{1}{8}}.
\end{equation}

\medskip

\underline{Step 3) Behavior of $\text{Tr}(\hat{A}^2)$:}

\medskip

In this step we shall use the celebrated Nourdin-Peccati bound which is at the heart of the Malliavin-Stein method. It asserts that, as soon as $X$ belongs to some Wiener chaos of any order and has variance $1$ we have

\begin{equation}\label{MallSteinBound}
d_{TV}\left(X,\mathcal{N}(0,1)\right)\le \frac{1}{\sqrt{3}} \sqrt{\E(X^4)-3}.
\end{equation}
Let us apply this bound to $\s F$, conditionally to $(\hat{X}_1,\cdots,\hat{X}_N)$. Taking into account that $\E(\s F^2)=2 \text{Tr}(\hat{A}^2)$, the bound (\ref{MallSteinBound}) gives, conditionally to $\hat{\mathcal{F}}$ that

\begin{eqnarray*}
\left|\E\left(e^{i\xi \s F}\right)-e^{-\xi^2 \text{Tr}(\hat{A}^2)}\right|&\le&d_{TV}\left(\s F, \mathcal{N}(0,2\text{Tr}(\hat{A}^2)\right)\\
&\le& \frac{1}{\sqrt{3}\text{Tr}(\hat{A}^2)}\sqrt{\kappa_4(\s F)}\\
&=&\frac{1}{\sqrt{3}\text{Tr}(\hat{A}^2)} \sqrt{48 \sum_{k=1}^N \hat{\lambda}_k^4}\\
&\le& 4 \frac{|\hat{\lambda}_1|}{\sqrt{\text{Tr}(\hat{A}^2)}}\end{eqnarray*}
Note that $\text{Tr}(\hat{A}^2)$ is a quadratic form in Gaussian random variables which satisfies

$$\hat{\E}\left(\text{Tr}(\hat{A}^2)\right)= 9 \sum_{1\le i_1,i_2,i_3\le N} a(i_1,i_2,i_3)^2 =\frac{3}{2}.$$

This is why we have $\hat{\mathbb{P}}\left(\text{Tr}(\hat{A}^2)=0\right)=0$ and the above inequalities are unambiguous. Taking the expectation with respect to $\hat{\E}$ gives

\begin{eqnarray*}
&&\left|\E\left(e^{-\xi^2 \frac{\Gamma[F,F]}{2}}\right)-\hat{\E}\left(e^{-\xi^2 \text{Tr}(\hat{A}^2)}\right)\right|\\
&=&\left|\hat{\E} \E\left(e^{i\xi \s F}\right) -\hat{\E}\left(e^{-\xi^2 \text{Tr}(\hat{A}^2)}\right)\right|\\
&\le&\hat{\E}\left(\left|\E\left(e^{i\xi \s F}\right)-e^{-\xi^2 \text{Tr}(\hat{A}^2)}\right|\right)\\
&\le& 4\hat{\E}\left(\frac{|\hat{\lambda}_1|}{\sqrt{\text{Tr}(\hat{A}^2)}}\right)\\
&\le& 4 \left(\hat{\E}\left(\frac{|\hat{\lambda}_1|}{\sqrt{\alpha}}\right)+\mathbb{P}\left(\text{Tr}(\hat{A}^2)\le \alpha\right)\right)
\end{eqnarray*}

However $\text{Tr}(\hat{A}^2)$ is polynomial of degree $2$ in the Gaussian r.v $(\hat{X}_1,\cdots,\hat{X}_N)$ and the Carbery-Wright inequality (see  e.g.\cite{carbery2001distributional}) gives, for some absolute constant $C>0$ which, in the sequel, may change from line to line that

$$\mathbb{P}\left(\text{Tr}(\hat{A}^2)\le \alpha\right) \le C \frac{\sqrt{\alpha}}{\hat{\E}\left(\text{Tr}(\hat{A}^2)\right)}\le c \sqrt{\alpha}.$$

Gathering all theses facts  leads to the inequalities valid for every $\xi\in\R$:

\begin{eqnarray*}
\left|\E\left(e^{-\xi^2 \frac{\Gamma[F,F]}{2}}\right)-\hat{\E}\left(e^{-\xi^2 \text{Tr}(\hat{A}^2)}\right)\right|&\le& C\left(\frac{1}{\sqrt{\alpha}}\hat{\E}\left(|\hat{\lambda}_1|\right)+\sqrt{\alpha}\right)\\
&\le& C \sqrt{\hat{\E}\left(|\hat{\lambda}_1|\right)} \\
&\le& C \kappa_4(F)^{\frac{1}{16}},
\end{eqnarray*}
where we have used the bound (\ref{finalboundspec}) at the end.

\medskip

Let us now show by the absurd that 
\begin{equation}\label{Gamma-dual}
\forall \epsilon>0,\exists \delta_\epsilon>0\,\,\text{s.t.}\,\,\kappa_4(F)<\delta_\epsilon \Rightarrow \text{Var}\left(\text{Tr}\left(\hat{A}^2\right)\right) \le \epsilon.
\end{equation}

If it would not be true, for some $\epsilon>0$ one might find a sequence $\{F_n\}_{n\ge 1}$ of random variables in the third chaos such that :

\begin{itemize}
\item $\E(F_n^2)=1$, 
\item $\kappa_4(F_n)\le \frac{1}{n}$
\item $\text{Var}\left(\text{Tr}\left(\hat{A_n}^2\right)\right) \ge \epsilon$. 
\end{itemize}

Up to extracting a subsequence, we may assume that $\text{Tr}\left(\hat{A_n}^2\right)$ converges in distribution towards $Z_\infty$ and one deduces by passing at the limit that, for every $\xi\in\R$:
$$\E\left(e^{-\xi^2 \frac{3}{2}}\right)=\hat{\E}\left(e^{-\xi^2  Z_\infty}\right).$$

Finally, by injectivity of the Laplace transform of positive random variables one deduces that $Z_\infty=\frac{3}{2}$ and thus necessarily $\text{Tr}\left(\hat{A}_n^2\right)\to \frac 3 2$ and we finally get
$$\text{Var}\left(\text{Tr}\left(\hat{A}_n^2\right)\right)\to 0,$$
which is contradictory. Unfortunately reasoning by the absurd leads to abstract constants and we lose here the ability of providing explicit quantitative statements. One way to bypass this problem is to relate the infinite norm of the difference of the Laplace transforms of $\Gamma[F,F]$ and $\text{Tr}\left(\hat{A}^2\right)$ to some usual probability distance, such as Forter-Mourier. This falls beyond the scope of this article and the statement (\ref{Gamma-dual}) would be enough for us.

\medskip

\underline{Step 4): Improving Carbery-Wright rate for $\Gamma[F,F]$ by a factor $3$:}

\medskip

Set $0<\theta<\frac{3}{4}$ and chose an integer $p$ large enough such that 
$$2\theta<\frac{1}{2}+\frac{p^2-p}{p^2+2p-2}<\frac{3}{2}.$$ The equation (\ref{Gamma-spec}) gives the following inequality:

\begin{eqnarray*}
\left|\E\left(e^{-\frac{\Gamma[F,F]}{2} \xi^2}\right)\right|&=&\left|\hat{\E}\left(\frac{1}{\sqrt{\det\left(\I-2 i \xi \hat{A}\right)}}\right)\right|\\
&\le& \hat{\E}\left(\prod_{k=1}^N \frac{1}{\left(1+4\hat{\lambda}_k^2 \xi^2\right)^{\frac 1 4}}\right).
\end{eqnarray*}

Let us use again some notations of the last section and set for every $q\ge 1$:

\begin{eqnarray*}
\hat{S}_{q}&=&\sum_{1\le i_1<i_2<\cdots<i_q\le N} \hat{\lambda}^2_{i_1}\hat{\lambda}^2_{i_2}\cdots \hat{\lambda}^2_{i_q}\\
\hat{\mathcal{N}}_q&=&\sum_{i=1}^N \hat{\lambda}_i^{2q}.
\end{eqnarray*}

As previously, the Newton-Girard inequalities give

$$\hat{S}_p=\frac{\hat{\mathcal{N}}_1^p}{p!}+(-1)^p \sum_{\begin{subarray}{c}
m_1<p\\m_1+2m_2+\cdots+pm_p=p
\end{subarray}} \prod_{i=1}^p \frac{(-\hat{\mathcal{N}}_i)^{m_i}}{m_i! i^{m_i}},$$

with $p$ being the integer previously chosen. Besides we also have for every choice of $(m_1,\cdots,m_p)$ such that $m_1+2m_2+\cdots+p m_p=p$ and $m_1<p$ there exits $i\ge 2$ such that $m_i>0$. Hence, for this index $i$ we have $|\hat{\mathcal{N}_i}|\le \left(|\hat{\lambda}_i|^{2i-2} \hat{\mathcal{N}}_1\right)$ which implies that

$$\prod_{i=1}^p |\hat{\mathcal{N}}_i|^{m_i}\le |\hat{\lambda}_1|^{\sum_{i=2}^N (2i-2) m_i} \hat{\mathcal{N}}_1^{m_1+\cdots+m_p}$$

Let us recall that one has $\hat{\E}\left(\hat{\mathcal{N}}_1\right)=\text{Tr}\left(\hat{A}^2\right)=\frac{3}{2}$ as well as the estimate (\ref{finalboundspec}). One deduces that, provided that $\kappa_4(F)$ is small enough, then

$$\hat{\E}\left(\hat{{S}}_p\right)\ge \frac{1}{2 p!} \hat{\E}\left(\hat{\mathcal{N}}_1^p\right)\ge \frac{1}{2} \frac{3^p}{2^p p!}.$$

On the other hand, the Newton-Girard formulas expresses $\hat{S}_p$ as linear combinations of $\text{Tr}(\hat{A}^k)$ with $k\in\{2,4,\cdots,2p\}$ which ensures that $\hat{S}_p$ is a polynomial function of degree $2p$ of the Gaussian variables $\{\hat{X}_1,\cdots,\hat{X}_N\}$ and which is positive. One can use the Carbery-Wright inequality to $\hat{S}_p$ which asserts that

\begin{equation}\label{CWSp}
\E\left(\hat{S}_p\right)^{\frac{1}{2p}}\hat{\mathbb{P}}\left(\hat{S}_p\le \alpha\right)\le C \alpha^{\frac{1}{2p}}.
\end{equation}

As a matter fact, we have shown that, provided that $\kappa_4(F)$ is small enough, one gets that

$$\mathbb{P}\left(\hat{S}_p \le \alpha\right) \le C_p \alpha^{\frac{1}{2p}},$$
where the constant $c_p$ is an absolute constant depending only on $p$. We are now ready for the final argument. Let us firs note that, since for every $p\ge 1$, $\hat{S}_p\ge 0$, the following inequality holds:
$$\prod_{k=1}^N \frac{1}{\left(1+4\hat{\lambda}_k^2 \xi^2\right)}\ge 1+ 4 \xi^2 \hat{S}_1+4^p \xi^{2p} \hat{S}_{p}.$$

Discussing according to $\hat{S}_{p}\ge \epsilon$ leads to
\begin{equation}\label{decoupage}
\E\left(\prod_{k=1}^N \frac{1}{\left(1+4\hat{\lambda}_k^2 \xi^2\right)^{\frac 1 4}}\right)\le \left(\frac{1}{4^p \xi^{2p} \epsilon}\right)^{\frac{1}{4}}+\frac{1}{\sqrt{\xi}}\E\left(\frac{1}{\hat{S}_1^{\frac 1 4}} \textbf{1}_{\{\hat{S}_{p}<\epsilon}\}\right)
\end{equation}

One the other hand, relying on the statement (\ref{Gamma-dual}) established in the previous step, we know that $\text{Var}(\text{Tr}(\hat{A}^2))$ tends to zero when $\kappa_4(F)$ tends to zero. Besides,

\begin{eqnarray*}
\text{Tr}(\hat{A}^2)=9 \sum_{1\le i_1,i_2\le N} \left(\sum_{k=1}^N a(i_1,i_2,k) \hat{X}_k\right)^2
\end{eqnarray*}

is a positive quadratic Gaussian form which can be diagonalized and written in the form

\begin{eqnarray*}
\text{Tr}(\hat{A}^2)=\sum_{k=1}^N \beta_k \hat{G}_k^2
\end{eqnarray*}
where $\{\beta_k\}_{k\ge 1}$ is a collection of positive real numbers whose sum is necessarily equal to $\hat{\E}\left(\text{Tr}\left(\hat{A}^2\right)\right)=\frac 3 2$ and where the $\hat{G}_k$ are independent standard Gaussian in the space $\left(\hat{\Omega},\mathcal{F},\mathbb{P}\right)$. Combining Theorem \ref{Mainsecond} with the observation (\ref{Gamma-dual}) entails that:
\begin{equation}\label{Integ-gammadual}
\forall q\ge 1,\exists \eta_q>0\,\,\text{s.t.},\,\,\kappa_4(F)\le \eta \Rightarrow \text{Tr}(\hat{A}^2)=\mathcal{N}_2=\hat{S}_1 \in L^{-q}(\hat{\Omega},\mathcal{F},\mathbb{P}).
\end{equation}
Assuming that $\kappa_4(F)<\eta_p$ and applying H\"{o}lder inequality gives, for a constant $C_p$ only depending on $p$ and changing from line to line:
\begin{eqnarray*}
\E\left(\prod_{k=1}^N \frac{1}{\left(1+4\hat{\lambda}_k^2 \xi^2\right)^{\frac 1 4}}\right)&\le& C_p \left(\frac{1}{\xi^{\frac{p}{2}} \epsilon^{\frac{1}{4}}}+\frac{1}{\sqrt{\xi}}\E\left(\frac{1}{\hat{S}_2^{\frac p 4}}\right)^{\frac{1}{p}}\hat{\mathbb{P}}\left(\hat{S}_{p}<\epsilon\right)^{\frac{p-1}{p}}\right)\\
&\le& C_p \left(\frac{1}{\xi^{\frac{p}{2}} \epsilon^{\frac{1}{4}}}+\frac{1}{\sqrt{\xi}}\left(\epsilon^{\frac{1}{2p}}\right)^{\frac{p-1}{p}}\right)
\end{eqnarray*}
When minimizing on $\epsilon$ one must have
$$\epsilon^{\frac{p-1}{p}\frac{1}{2p}+\frac{1}{4}}=\frac{1}{\xi^{\frac{p-1}{2}}},$$
which gives
$$\epsilon^{\frac{1}{2p}}=\left(\frac{1}{\xi}\right)^{\frac{p^2-p}{p^2+2p-2}}.$$
Hence, one gets for $\xi>0$ large enough:
$$
\E\left(\prod_{k=1}^N \frac{1}{\left(1+4\hat{\lambda}_k^2 \xi^2\right)^{\frac 1 4}}\right)\le C_p \left(\frac{1}{\xi}\right)^{\frac{p^2-p}{p^2+2p-2}+\frac 1 2}\le C_p \left(\frac{1}{\xi}\right)^{2\theta}.
$$
To conclude the proof we use again Markov inequality and write
\begin{eqnarray*}
\mathbb{P}\left(\Gamma[F,F]<\epsilon\right)&=&\mathbb{P}\left(e^{-\xi^2 \Gamma[F,F]}\ge e^{-\xi^2\epsilon}\right)\\
&\le& e^{\xi^2 \epsilon} \E\left(e^{-\xi^2 \Gamma[F,F]}\right)\\
&\le& C_p e^{\xi^2 \epsilon} \frac{1}{\xi^{2\theta}}.
\end{eqnarray*}
Choosing $\xi=\frac{1}{\sqrt{\epsilon}}$ leads to the desired conclusion.
\end{proof}
\bibliographystyle{alpha}
\bibliography{ref-fmt}

\end{document}